\documentclass{ifcolog-c}
\usepackage{bussproofs}
\usepackage{mathrsfs}
\usepackage{amsthm}
\usepackage{amssymb}
\usepackage[dvipsnames]{xcolor}
\usepackage{stmaryrd}
\usepackage{calc}
\usepackage{cmll}
\usepackage{amsmath}
\usepackage{array}
\newtheorem{theorem}{Theorem}
\newtheorem{definition}{Definition}
\newtheorem{proposition}{Proposition}

\title{A note on schematic validity and completeness in Prawitz's semantics}
\titlethanks{I thank Cesare Cozzo, Ansten Klev, Ivo Pezlar, Thomas Piecha, Paolo Pistone, Dag Prawitz, Peter Schroeder-Heister and Will Stafford for helpful comments. Work on this article was supported by grant LQ300092101 from the Czech Academy of Sciences}
\addauthor[piccolomini@flu.cas.cz]{Antonio Piccolomini d'Aragona}{Institute of Philosophy, Czech Academy of Sciences}
\authorrunning{Piccolomini d'Aragona}
\titlerunning{A note on schematicity and completeness in Prawitz}

\date{March 2021}

\begin{document}

\nopagenumber

\maketitle

\noindent\small \textbf{Abstract.} 
I discuss two approaches to monotonic proof-theoretic semantics. In the first one, which I call SVA, consequence is understood in terms of existence of valid arguments. The latter involve the notions of argument structure and justification for arbitrary non-introduction rules. In the second approach, which I call Base Semantics, structures and justifications are left aside, and consequence is defined outright over background atomic theories. Many (in)completeness results have been proved relative to Base Semantics, the question being whether these can be extended to SVA. By limiting myself to a framework with classical meta-logic, I prove correctness of classical logic on Base Semantics, and show that this result adapts to SVA when justifications are allowed to be choice-functions over atomic theories or unrestricted reduction systems of argument structures. I also point out that, however, if justifications are required to be more schematic, correctness of classical logic over SVA may fail, even with classical logic in the meta-language. This seems to reveal that the way justifications are understood may be a distinguishing feature of different accounts of proof-theoretic validity.


\section{Introduction}

Prawitz's semantics, an instance of proof-theoretic semantics \cite{francez, schroeder-heisterSE}, has come in basically two forms: \emph{semantics of valid arguments} -- SVA, see e.g. \cite{prawitz73} -- and \emph{theory of grounds} -- ToG, see e.g. \cite{prawitz2015}. Here, I shall focus mostly on SVA, and I will only occasionally refer to ToG.

SVA is based on the notion of \emph{valid argument}. The latter is inspired by Prawitz's normalisation results for Gentzen's Natural Deduction \cite{prawitz65}, stating that derivations for $\Gamma \vdash A$ can be transformed, through suitable reductions, to derivations for $\Gamma^* \subseteq \Gamma \vdash A$ without \emph{detours}. A detour is given by a formula which occurs both as conclusion of an introduction, and as a major premise of an elimination.

In intuitionistic logic, Prawitz's theorems imply what Schroeder-Heister called the \emph{fundamental corollary} \cite{schroeder-heister2006}: $A$ is a theorem iff there is a closed derivation of $A$ ending by an introduction. This may confirm Gentzen's claim that introductions fix meaning, while eliminations are unique functions of the introductions \cite{gentzen}.

Prawitz's well-known \emph{inversion principle} \cite{prawitz65} that a by-introduction proof of $A$ contains (part of) what is needed for drawing consequences of $A$, can thus undergo a semantic reading: arbitrary inferences may be \emph{justified} by transformations on proofs of the major premises. An argument in general can be said to be valid when its arbitrary inferences are so justified. This might happen relative to some underlying atomic proof-system, meant to determine non-logical meanings, or hold irrespective of such systems, and thus qualify as logical.

In the SVA-inspired framework of some other authors, which following Sandqvist \cite{sandqvistcompl} we may call \emph{Base Semantics}, the role of justifications and proof-structures is limited to atomic proof-systems. Validity of arguments is replaced by a consequence relation among formulas, which may once again be relative to given atomic proof-systems, or hold on all such systems, i.e. be logical.

Works in the Base Semantics tradition have shown that intuitionistic completeness depends on peculiar features of atomic proof-systems -- see \cite{piecha} for an overview, while later examples are \cite{piechaschroeder-heister, schroeder-heisterolf, stafford1, stafford2}. Building upon a result of \cite{piechaschroeder-heister}, I here prove incompleteness of intuitionistic logic with respect to a variant of SVA, and so show that, under a certain reading, Base Semantics copes with an approach where justifications and proof-structures are not disregarded. However, I also point out that the passage from Base Semantics to SVA is non-trivial, since the adaptation of the result of \cite{piechaschroeder-heister} to SVA seems to force an understanding of reductions as \emph{non-schematic} functions from and to proof-structures.

In particular, I shall argue that a more ‘‘schematic'' understanding of reductions blocks certain proofs of soundness of classical logic, which instead obtain when reductions are non-schematic. This may mean that a substantive part of the constructive nature of Prawitz's semantics stems from a certain understanding of reductions, rather than just from the fact that the approach is proof-based.

Although SVA and Base Semantics can be (and have been) developed for full first-order logic, intuitionistic (in)completeness is mostly discussed at a propositional level. I accordingly limit myself to propositional logic. Moreover, in proving intuitionistic incompleteness, I use classical logic at the meta-level. This is however sufficient for raising my point.

The article is structured as follows. In Section 2 I give an outline of (a variant of) SVA. In Section 3 I define (a variant of) Base Semantics, and prove incompleteness of intuitionistic logic with respect to it -- specifically, soundness of excluded middle with respect to it. In Section 4 I adapt the incompleteness proof to SVA in such a way as to stick to Prawitz's original idea of a semantics where justifications and proof-structures do play an active role. In Section 5 I discuss the incompleteness result of Section 4 with reference to a notion of justification understood on different degrees of strength. In the concluding remarks, I finally suggest some tentative proposals for accounting for a ‘‘schematic'' notion of reduction, in terms of restrictions on the form of open valid arguments, of Pezlar selector, and of linearity conditions on replacement of variable parts in proof-structures. 

\section{SVA (over a base)}

Prawitz's normalisation obtains through reduction functions from and to Gentzen's Natural Deduction derivations, whose iteration eventually normalise application arguments while keeping the conclusion and not expanding the assumptions-set. For disjunction detours, such a function is e.g.
\begin{prooftree}
    \AxiomC{$\mathscr{D}_1$}
    \noLine
    \UnaryInfC{$A_i$}
    \RightLabel{($\vee_I$), $i = 1, 2$}
    \UnaryInfC{$A_1 \vee A_2$}
    \AxiomC{$[A_1]$}
    \noLine
    \UnaryInfC{$\mathscr{D}_2$}
    \noLine
    \UnaryInfC{$B$}
    \AxiomC{$[A_2]$}
    \noLine
    \UnaryInfC{$\mathscr{D}_3$}
    \noLine
    \UnaryInfC{$B$}
    \RightLabel{($\vee_E$)}
    \TrinaryInfC{$B$}
    \AxiomC{$\stackrel{\phi_\vee}{\Longrightarrow}$}
    \noLine
    \UnaryInfC{}
    \AxiomC{$\mathscr{D}_1$}
    \noLine
    \UnaryInfC{$[A_i]$}
    \noLine
    \UnaryInfC{$\mathscr{D}_{i + 1}$}
    \noLine
    \UnaryInfC{$B$}
    \noLine
    \TrinaryInfC{}
\end{prooftree}
This is brought to a semantic level by SVA, through a suitable generalisation of the notions of derivation and reduction towards, respectively, the notions of argument structure and justification, relative to a background language $\mathscr{L}$ which will be here given as follows.
\begin{definition}
The grammar of the \emph{language} $\mathscr{L}$ is 
\begin{center}
    $X \coloneqq p \ | \ \bot \ | \ X \wedge X \ | \ X \vee X \ | \ X \rightarrow X$.
\end{center}
\end{definition}
\noindent The $p$-s and $\bot$ are atoms, $\bot$ being a constant symbol for the absurd. As done above, I use capital Latin letter as variables for formulas of $\mathscr{L}$, while capital Greek letters will be used as variables for sets of formulas. Negation is defined by
\begin{center}
    $\neg A \stackrel{def}{=} A \rightarrow \bot$.
\end{center}

Just like in Model Theory, also in SVA we have structures which fix the meaning of the non-logical components of $\mathscr{L}$, and thus serve as induction-base for the definition of the semantic core-notions. Unlike Model Theory, however, these structures are not sets (or values) onto which non-logical components are mapped, but atomic systems, i.e. sets of rules that establish the semantic behaviour of atoms -- and thus determine constructively the meaning of the components they involve. Atomic systems can be defined in different ways, depending on different desiderata one may have -- see \cite{piechaschroeder-heisterbasis} -- and yielding different outcomes for validity-related concepts such as completeness -- see \cite{piecha}. Here, I will understand them as sets of \emph{production rules}, i.e. I will not allow premises to be lower-level rules, nor will I allow for discharge of assumptions at the atomic level.
\begin{definition}
An \emph{atomic base} $\mathfrak{B}$ over $\mathscr{L}$ is a (countable) set of rules
\begin{prooftree}
    \AxiomC{$A_1, ..., A_n$}
    \UnaryInfC{$B$}
\end{prooftree}
where $n \geq 0$, $A_i, B$ are atoms of $\mathscr{L}$ and $A_i \neq \bot$ ($i \leq n$).
\end{definition}
\noindent Derivations on $\mathfrak{B}$ are defined in a usual inductive way. The derivability relation on $\mathfrak{B}$ is written $\vdash_{\mathfrak{B}}$, while the derivations-set of $\mathfrak{B}$ is written $\texttt{DER}_{\mathfrak{B}}$. I always require bases to be consistent.
\begin{definition}
    $\mathfrak{B}$ is \emph{consistent} iff $\not\vdash_{\mathfrak{B}} \bot$.
\end{definition}
\begin{definition}
An \emph{argument structure} over $\mathscr{L}$ is a pair $\langle T, f \rangle$ such that $T$ is a tree whose nodes are either empty, in which case they are always top-nodes, or formulas of $\mathscr{L}$, and $f$ is a function defined on a subset $\Gamma$ of non-empty top-nodes of $T$ such that, for every $A \in \Gamma$, $f(A)$ is below $A$ in $T$.
\end{definition}
\noindent Given $\mathscr{D} = \langle T, f \rangle$ where $T$ has top-nodes $\Gamma$ and root $A$, I call $\Gamma$ the \emph{assumptions} of $\mathscr{D}$ and $A$ the \emph{conclusion} of $\mathscr{D}$ -- thus $f$ is an assumptions-discharge.
\begin{definition}
    $\mathscr{D}$ is \emph{closed} iff all its assumptions are discharged, otherwise it is \emph{open}.
\end{definition}
\noindent Where $\Gamma$ is the set of the open assumptions of a $\mathscr{D}$ with conclusion $A$, I shall say that $\mathscr{D}$ is an argument structure \emph{from} $\Gamma$ \emph{to} $A$, and I shall indicate it by
\begin{prooftree}
    \AxiomC{$\Gamma$}
    \noLine
    \UnaryInfC{$\mathscr{D}$}
    \noLine
    \UnaryInfC{$A$}
\end{prooftree}

\begin{definition}
    Given $\mathscr{D}$ from $\Gamma$ to $A$ and $\sigma$ a function such that, for every $B \in \Gamma$, $\sigma(B)$ is a (closed) argument structure with conclusion $B$, a \emph{(closed) $\sigma$-instance} $\mathscr{D}^\sigma$ of $\mathscr{D}$ is the argument structure obtained from $\mathscr{D}$ by replacing every $B \in \Gamma$ with $\sigma(B)$.
\end{definition}

\begin{definition}
    An \emph{inference} is a triple $\langle \langle \mathscr{D}_1, ..., \mathscr{D}_n \rangle, A, \delta \rangle$, where $\delta$ is an indication of assumptions which may be discharged by the inference. The argument structure \emph{associated} to the inference, indicated by
    \begin{prooftree}
        \AxiomC{$\mathscr{D}_1, ..., \mathscr{D}_n$}
        \RightLabel{$\delta$}
        \UnaryInfC{$A$}
    \end{prooftree}
    is obtained by conjoining the trees of the $\mathscr{D}_i$-s through a root-node $A$, and by expanding the assumptions-discharges of the $\mathscr{D}_i$-s according to $\delta$. A \emph{rule} is a set of inferences, whose elements are called \emph{instances} of the rule.
\end{definition}
\noindent I shall assume that rules can be described schematically, e.g. standard introduction rules in Gentzen's Natural Deduction,
\begin{prooftree}
    \AxiomC{$A$}
    \AxiomC{$B$}
    \RightLabel{($\wedge_I$)}
    \BinaryInfC{$A \wedge B$}
    \AxiomC{$A_i$}
    \RightLabel{($\vee_I$), $i = 1, 2$}
    \UnaryInfC{$A_1 \vee A_2$}
    \AxiomC{$[A]$}
    \noLine
    \UnaryInfC{$B$}
    \RightLabel{($\rightarrow_I$)}
    \UnaryInfC{$A \rightarrow B$}
    \noLine
    \TrinaryInfC{}
\end{prooftree}

\begin{definition}
    $\mathscr{D}$ is \emph{canonical} iff it is associated to an instance of an introduction rule, otherwise it is \emph{non-canonical}.
\end{definition}

\begin{definition}
    Given a rule $R$, a \emph{justification} of $R$ is a constructive function $\phi$ defined on the set of the argument structures $\mathbb{D}$ associated to some sub-set of $R$ such that, for every $\mathscr{D} \in \mathbb{D}$,
    \begin{itemize}
        \item $\mathscr{D}$ is from $\Gamma$ to $A \Longrightarrow \phi(\mathscr{D})$ is from $\Gamma^* \subseteq \Gamma$ to $A$, and
        \item for every $\sigma$, $\phi$ is defined on $\mathscr{D}^\sigma$ and $\phi(\mathscr{D}^\sigma) = \phi(\mathscr{D})^\sigma$.
    \end{itemize}
\end{definition}
\noindent For example, the reduction $\phi_\vee$ for eliminating maximal disjunctions from Natural Deduction derivations, which has been defined at the beginning of this section, can be understood semantically as a justification showing that elimination of disjunction ($\vee_E$) can be safely removed, \emph{salva} provability, when its major premise has been obtained by ($\vee_I$). Observe that $\phi_\vee$ needs not be defined on all the applications of ($\vee_E$). Via Definition 12 below, applications whose major premise is proved via ($\vee_I$) only suffice for showing the rule to be justified in the sense hinted at above.

I take the notion of (\emph{immediate}) \emph{sub-structure} of $\mathscr{D}$ to be clear enough at this point, so I will not define it explicitly. The same I will do with the notion of \emph{substitution} of a sub-structure $\mathscr{D}^*$ by an argument structure $\mathscr{D}^{**}$ in an argument structure $\mathscr{D}$, written $\mathscr{D}[\mathscr{D}^{**}/\mathscr{D}^*]$ -- observe that this already occurs in Definition 6. The application of substitution may require re-indexing the discharge functions associated to the argument structures, but I shall not deal with these details here.
\begin{definition}
    Given a set of justifications $\mathfrak{J}$, $\mathscr{D}$ \emph{immediately reduces} to $\mathscr{D}^*$ \emph{relative to} $\mathfrak{J}$, written $\mathscr{D} \leq^{\mathfrak{J}}_\iota \mathscr{D}^*$, iff $\mathscr{D} = \mathscr{D}^*$ or, for some sub-structure $\mathscr{D}^{**}$ of $\mathscr{D}$ and some $\phi \in \mathfrak{J}$, $\phi$ is defined on $\mathscr{D}^{**}$ and $\mathscr{D}^* = \mathscr{D}[\phi(\mathscr{D}^{**})/\mathscr{D}^{**}]$. $\mathscr{D}$ \emph{reduces} to $\mathscr{D}^*$ \emph{relative to} $\mathfrak{J}$, written $\mathscr{D} \leq^{\mathfrak{J}} \mathscr{D}^*$, iff there is a sequence $\mathscr{D} = \mathscr{D}_1 \leq^{\mathfrak{J}}_\iota \mathscr{D}_2 \leq^{\mathfrak{J}}_\iota ... \leq^{\mathfrak{J}}_\iota \mathscr{D}_{n - 1} \leq^\mathfrak{J}_\iota \mathscr{D}_n = \mathscr{D}^*$.
\end{definition}

\begin{definition}
    An \emph{argument} is a pair $\langle \mathscr{D}, \mathfrak{J} \rangle$.
\end{definition}

\begin{definition}
    $\langle \mathscr{D}, \mathfrak{J} \rangle$ is \emph{valid on} $\mathfrak{B}$ iff
    \begin{itemize}
    \item $\mathscr{D}$ is closed $\Longrightarrow$ 
    \begin{itemize}
        \item the conclusion of $\mathscr{D}$ is an atom $\Longrightarrow \mathscr{D} \leq^\mathfrak{J} \mathscr{D}^*$ for $\mathscr{D}^* = \langle T, \emptyset \rangle$ and $T \in \texttt{DER}_\mathfrak{B}$ closed;
        \item the conclusion of $\mathscr{D}$ is not an atom $\Longrightarrow \mathscr{D} \leq^\mathfrak{J} \mathscr{D}^*$ for $\mathscr{D}^*$ closed canonical with immediate sub-structures valid on $\mathfrak{B}$ when paired with $\mathfrak{J}$;
    \end{itemize}
        \item $\mathscr{D}$ is open from $\Gamma$ to $A \Longrightarrow$ for every $\sigma$, every $B \in \Gamma$ and every extension $\mathfrak{J}^+$ of $\mathfrak{J}$, if $\langle \sigma(B), \mathfrak{J}^+ \rangle$ is valid on $\mathfrak{B}$, then $\langle \mathscr{D}^\sigma, \mathfrak{J}^+ \rangle$ is valid on $\mathfrak{B}$.
    \end{itemize}
\end{definition}
\noindent Let us prove that $\langle \mathscr{D}, \{\phi_\vee\} \cup \mathfrak{J}_2 \cup \mathfrak{J}_3 \rangle$ is valid on any $\mathfrak{B}$, where $\mathscr{D}$ is the structure
\begin{prooftree}
    \AxiomC{$A_1 \vee A_2$}
    \AxiomC{$[A_1]$}
    \noLine
    \UnaryInfC{$\mathscr{D}_2$}
    \noLine
    \UnaryInfC{$B$}
    \AxiomC{$[A_2]$}
    \noLine
    \UnaryInfC{$\mathscr{D}_3$}
    \noLine
    \UnaryInfC{$B$}
    \RightLabel{($\vee_E$)}
    \TrinaryInfC{$B$}
\end{prooftree}
and $\langle \mathscr{D}_2, \mathfrak{J}_2 \rangle$ and $\langle \mathscr{D}_3, \mathfrak{J}_3 \rangle$ are open valid on $\mathfrak{B}$ from $A_1$ and $A_2$ respectively to $B$. By the second clause of Definition 12, we must prove that, for every $\langle \mathscr{D}_1, (\{\phi_\vee\} \cup \mathfrak{J}_2 \cup \mathfrak{J}_3)^+ \rangle$ valid on $\mathfrak{B}$ where $\mathscr{D}_1$ is closed and has conclusion $A_1 \vee A_2$, we have that $\langle \mathscr{D}[\mathscr{D}_1/A_1 \vee A_2], (\{\phi_\vee\} \cup \mathfrak{J}_2 \cup \mathfrak{J}_3)^+ \rangle$ is valid on $\mathfrak{B}$. Since $\langle \mathscr{D}_1, (\{\phi_\vee\} \cup \mathfrak{J}_2 \cup \mathfrak{J}_3)^+ \rangle$ is closed valid on $\mathfrak{B}$, by the first clause of Definition 12, there is a closed canonical $\mathscr{D}^*_1$ with conclusion $A_1 \vee A_2$ such that $\mathscr{D}_1 \leq_{(\{\phi_\vee\} \cup \mathfrak{J}_2 \cup \mathfrak{J}_3)^+} \mathscr{D}^*_1$, and $\langle \mathscr{D}^*_1, (\{\phi_\vee\} \cup \mathfrak{J}_2 \cup \mathfrak{J}_3)^+ \rangle$ is valid on $\mathfrak{B}$. Also, $\mathscr{D}_1 \leq_{(\{\phi_\vee\} \cup \mathfrak{J}_2 \cup \mathfrak{J}_3)^+} \mathscr{D}^*_1$ obviously implies $\mathscr{D} \leq_{(\{\phi_\vee\} \cup \mathfrak{J}_2 \cup \mathfrak{J}_3)^+} \mathscr{D}^*$, where $\mathscr{D}^*$ is closed of the form
\begin{prooftree}
    \AxiomC{$\mathscr{D}^*_1$}
    \noLine
    \UnaryInfC{$A_i$}
    \RightLabel{($\vee_I$), $i = 1, 2$}
    \UnaryInfC{$A_1 \vee A_2$}
    \AxiomC{$[A_1]$}
    \noLine
    \UnaryInfC{$\mathscr{D}_2$}
    \noLine
    \UnaryInfC{$B$}
    \AxiomC{$[A_2]$}
    \noLine
    \UnaryInfC{$\mathscr{D}_3$}
    \noLine
    \UnaryInfC{$B$}
    \RightLabel{($\vee_E$)}
    \TrinaryInfC{$B$}
\end{prooftree}
We must show that $\langle \mathscr{D}^*, (\{\phi_\vee\} \cup \mathfrak{J}_2 \cup \mathfrak{J}_3)^+ \rangle$ is valid on $\mathfrak{B}$, since this implies $\langle \mathscr{D}[\mathscr{D}_1/A_1 \vee A_2], (\{\phi_\vee\} \cup \mathfrak{J}_2 \cup \mathfrak{J}_3)^+ \rangle$ valid on $\mathfrak{B}$ which in turn, by the arbitrary choice of $\langle \mathscr{D}_1, (\{\phi_\vee\} \cup \mathfrak{J}_2 \cup \mathfrak{J}_3)^+ \rangle$, implies our result. Since $\phi_\vee \in (\{\phi_\vee\} \cup \mathfrak{J}_2 \cup \mathfrak{J}_3)^+$, we have that $\mathscr{D}^* \leq_{(\{\phi_\vee\} \cup \mathfrak{J}_2 \cup \mathfrak{J}_3)^+} \mathscr{D}^{**}$, with $\mathscr{D}^{**}$ closed of the form
\begin{prooftree}
    \AxiomC{$\mathscr{D}^*_1$}
    \noLine
    \UnaryInfC{$[A_i$]}
    \noLine
    \UnaryInfC{$\mathscr{D}_{i + 1}$}
    \noLine
    \UnaryInfC{$B$}
\end{prooftree}
Since we assumed $\langle \mathscr{D}_{i + 1}, \mathfrak{J}_{i + 1} \rangle$ to be open valid on $\mathfrak{B}$, by the second clause of Definition 12 $\langle \mathscr{D}^{**}, (\{\phi_\vee\} \cup \mathfrak{J}_2 \cup \mathfrak{J}_3)^+ \rangle$ is closed valid on $\mathfrak{B}$, whence we are done.

\section{Base semantics}

Since argumental validity essentially depends on atomic provability on the base, on reducibility of non-canonical configurations to suitable introduction forms, and on assumptions-closure in the open-arguments case, we may decide to prune Definition 12 by just focusing on formulas, i.e. by dropping argument structures and justifications out.

This leads to Base Semantics, which I give in a somewhat simplified version as what \cite{piechaschroeder-heister} calls \emph{non-extension} semantics.
\begin{definition}
    $\Gamma \models_{\mathfrak{B}} A$ iff:
    \begin{itemize}
        \item[\emph{(a)}] $\Gamma = \emptyset \Longrightarrow$
        \begin{itemize}
            \item[\emph{(At)}] $A$ is an atom $\Longrightarrow \ \vdash_{\mathfrak{B}} A$;
            \item[\emph{($\wedge$)}] $A = B \wedge C \Longrightarrow \ \models_{\mathfrak{B}} B$ and $\models_{\mathfrak{B}} C$;
            \item[\emph{($\vee$)}] $A = B \vee C \Longrightarrow \ \models_{\mathfrak{B}} B$ or $\models_{\mathfrak{B}} C$;
            \item[\emph{($\rightarrow$)}] $A = B \rightarrow C \Longrightarrow B \models_\mathfrak{B} C$;
        \end{itemize}
        \item[\emph{(b)}] $\Gamma \neq \emptyset \Longrightarrow (\models_\mathfrak{B} \Gamma \Longrightarrow \ \models_\mathfrak{B} A)$, where $\models_\mathfrak{B} \Gamma$ means $\models_\mathfrak{B} B$ for every $B \in \Gamma$.
    \end{itemize}
\end{definition}
\begin{definition}
    $\Gamma \models A$ iff, for every $\mathfrak{B}$, $\Gamma \models_\mathfrak{B} A$.
\end{definition}

\noindent We can now give an easy incompleteness proof for intuitionistic logic -- $\texttt{IL}$. In fact, we can even prove soundness of classical logic (which is expected by the use of the latter in the meta-language).

\begin{proposition}
    For every $\mathfrak{B}$, $\models_\mathfrak{B} A \vee \neg A$.
\end{proposition}
\begin{proof}
    For $\mathfrak{B}$ arbitrary, by using classical logic in the meta-language, either $\models_\mathfrak{B} A$ or $\not\models_\mathfrak{B} A$. If $\models_\mathfrak{B} A$, then $\models_\mathfrak{B} A \vee \neg A$ by ($\vee$). If instead $\not\models_\mathfrak{B} A$, then $A \models_\mathfrak{B} \bot$ holds vacuously by (b). Hence, $\models_\mathfrak{B} \neg A$ by ($\rightarrow$), and $\models_\mathfrak{B} A \vee \neg A$ by ($\vee$). The result now follows from the arbitrariness of $\mathfrak{B}$.
\end{proof}

\begin{proposition}
    $\models A \vee \neg A$.
\end{proposition}

\begin{theorem}
    There are $\Gamma$ and $A$ such that $\Gamma \models A$ and $\Gamma \not\vdash_{\texttt{IL}} A$.
\end{theorem}
\noindent Except for the reference to atomic bases -- which is in a sense peculiar to Prawitzian approaches -- this is essentially the same as what one obtains in BHK-semantics with classical meta-language.

The crucial results of intuitionistic incompleteness to be found in literature on Base Semantics are of course much more significant than the one provided by Theorem 1 -- and, accordingly, their proofs require a much more fine-grained framework than the one I put forward here. In \cite{piechaschroeder-heister}, incompleteness of intuitionistic logic is e.g. referred to \emph{Harrop rule}, whose validity via Definition 14 is proved without using classical logic in the meta-language.\footnote{It should be however remarked that, in \cite{piechaschroeder-heister}, the rule from premise $A \rightarrow (B \vee C)$ to conclusion $(A \rightarrow B) \vee (A \rightarrow C)$ -- with no restrictions on $A$ -- is proved to be valid according to Definition 14, by using classical meta-logic and (b) in Definition 13. This can be easily adapted to the SVA-framework along the lines of Section 4 below. I shall not dwell upon this here, as the use of classical meta-logic validates excluded middle, which immediately implies intuitionistic incompleteness (both in Base Semantics and in a suitably modified version of SVA).} The adaptation of the Base Semantics framework to SVA that I present in Section 4 below can be extended to these more refined approaches; likewise, the critical remarks about the principles of Base Semantics I provide in Section 5 below also holds for more detailed accounts -- both in Base Semantics and in adaptations of it to SVA. Thus, the toy-example that, in the (trivialising) context of a classical meta-language, I discuss in Sections 4 and 5, are enough for raising my points. A comprehensive treatment of the issue can be developed in future works.\footnote{Let me just remark that the fundamental results proved in \cite{piechaschroeder-heister} require some additional principles, called \emph{Import} and \emph{Export}, and possibly higher-level rules at the atomic level -- i.e. atomic rules where the discharge of assumptions is allowed. The \emph{Import} principle says that, roughly, some assumptions may generate new atomic rules: $\Gamma \models_\mathfrak{B} A \Leftrightarrow \ \models_{\mathfrak{B} \cup \Sigma^\Gamma} A$. Together with higher-level atomic rules, this yields validity of Harrop rule via Definition 14 in non-extensions semantics. The \emph{Export} principle says instead that, roughly, (first-level) atomic rules generate new assumptions: $\Gamma \models_{\mathfrak{B} \cup \Sigma} A \Leftrightarrow \Gamma, \Delta^\Sigma \models_\mathfrak{B} A$. This yields validity of Harrop rule via Definition 14 for \emph{extension} semantics, namely Base Semantics where condition (b) in Definition 13 is given a monotonic form: $\Gamma \models_\mathfrak{B} A \Leftrightarrow$ for all $\mathfrak{C} \supseteq \mathfrak{B}, (\models_\mathfrak{C} \Gamma \Rightarrow \ \models_\mathfrak{C} A)$.}

\section{A kind of intuitionistic SVA-incompleteness}

The obvious strategy -- also suggested by \cite{piechaschroeder-heister} -- for moving from Base Semantics to SVA, is that of starting from the following equivalence:
\begin{itemize}
    \item[(EQ)] $\Gamma \models_\mathfrak{B} A$ iff there is $\langle \mathscr{D}, \mathfrak{J} \rangle$ from $\Gamma$ to $A$ valid on $\mathfrak{B}$.
\end{itemize}
Condition (a) in Definition 13 holds in SVA under (EQ), so we may just translate Proposition 1 via (EQ), after observing that, with classical logic in the meta-language, condition (b) in Definition 13 also holds under (EQ) in SVA -- proof omitted.

However, I prefer to stick to a strict SVA-formulation, which will allow me to remark some points that might be concealed when argument structures and justifications are entirely dropped out.

\begin{proposition}
    For every $\mathfrak{B}$, there is a closed argument for $A \vee \neg A$ valid on $\mathfrak{B}$.
\end{proposition}
\begin{proof}
    Let $\mathfrak{B}$ be arbitrary. By reasoning classically, either there is a closed valid argument for $A$, or there is not. Suppose there is not, and consider
    \begin{prooftree}
        \AxiomC{$A$}
        \LeftLabel{$\mathscr{D}^* =$}
        \UnaryInfC{$\bot$}
    \end{prooftree}
    Then, $\langle \mathscr{D}^*, \emptyset \rangle$ is vacuously valid on $\mathfrak{B}$, and so is $\langle \mathscr{D}^{**},\emptyset \rangle$, where $\mathscr{D}^{**}$ is
    \begin{prooftree}
    \AxiomC{$1$}
    \noLine
    \UnaryInfC{$[A]$}
    \UnaryInfC{$\bot$}
    \RightLabel{($\rightarrow_I$), $1$}
    \UnaryInfC{$A \rightarrow \bot$}
    \end{prooftree}
    Then, $\langle \mathscr{D}, \{\kappa_1\} \rangle$ is valid on $\mathfrak{B}$, where $\mathscr{D}$ is
    \begin{prooftree}
        \AxiomC{}
        \UnaryInfC{$A \vee \neg A$}
    \end{prooftree}
    and $\kappa_1$ is the constant function defined by
    \begin{prooftree}
        \AxiomC{$\mathscr{D}$}
        \noLine
        \UnaryInfC{}
        \AxiomC{$\Longrightarrow$}
        \noLine
        \UnaryInfC{}
        \AxiomC{$\mathscr{D}^{**}$}
        \noLine
        \UnaryInfC{$\neg A$}
        \RightLabel{($\vee_I$)}
        \UnaryInfC{$A \vee \neg A$}
        \noLine
        \TrinaryInfC{}
    \end{prooftree}
    Suppose there is a closed $\langle \mathscr{D}^{***}, \mathfrak{J} \rangle$ for $A$ valid on $\mathfrak{B}$. Then $\langle \mathscr{D}, \{\kappa_2\} \cup \mathfrak{J} \rangle$ is valid on $\mathfrak{B}$, where $\kappa_2$ is the constant function defined by
    \begin{prooftree}
        \AxiomC{$\mathscr{D}$}
        \noLine
        \UnaryInfC{}
        \AxiomC{$\Longrightarrow$}
        \noLine
        \UnaryInfC{}
        \AxiomC{$\mathscr{D}^{***}$}
        \noLine
        \UnaryInfC{$A$}
        \RightLabel{($\vee_I$)}
        \UnaryInfC{$A \vee \neg A$}
        \noLine
        \TrinaryInfC{}
    \end{prooftree}
    The result now follows by arbitrariness of $\mathfrak{B}$.
\end{proof}
\noindent A similar result for classical \emph{reduction ad absurdum} can be found in \cite{prawitz73}. Observe that, in Proposition 3, the role of (b) in Definition 13 is played by empty sets of justifications or by constant functions. Definition 14 via (EQ) now gives a kind incompleteness of $\texttt{IL}$.
\begin{definition}
    $\Gamma \models_{\Delta} A$ iff, for every $\mathfrak{B}$, there is $\langle \mathscr{D}, \mathfrak{J} \rangle$ from $\Gamma$ to $A$ valid on $\mathfrak{B}$.
\end{definition}
\begin{proposition}
    $\models_{\Delta} A \vee \neg A$.
\end{proposition}
\begin{theorem}
    There are $\Gamma$ and $A$ such that $\Gamma \models_{\Delta} A$ and $\Gamma \not\vdash_{\texttt{IL}} A$.
\end{theorem}

\section{Schematicity}

When read against the original formulation of SVA in \cite{prawitz73}, the proof of Proposition 3 implies a number of aspects which may be so modified as to obtain a different reading of $\models_\Delta$ than the one I have proposed in the previous section. First, we may invert the quantifiers in Definition 15.
\begin{definition}
    $\Gamma \models^*_\Delta A$ iff, for some $\langle \mathscr{D}, \mathfrak{J} \rangle$ from $\Gamma$ to $A$ and every $\mathfrak{B}$, $\langle \mathscr{D}, \mathfrak{J} \rangle$ is valid on $\mathfrak{B}$.
\end{definition}
\noindent Our proof for Theorem 2 does not directly apply to $\models^*_\Delta$. Yet, the proof of Proposition 3 shows that, for every $\mathfrak{B}$, there is a set of justifications $\mathfrak{J}^\mathfrak{B}$ such that $\langle \mathscr{D}, \mathfrak{J}^\mathfrak{B} \rangle$ valid on $\mathfrak{B}$, where $\mathscr{D}$ is the axiom $A \vee \neg A$, i.e.
\begin{center}
    $\mathfrak{J}^\mathfrak{B} = \begin{cases} \{\kappa_1\} & \text{if} \ \not\models_\mathfrak{B} A \\ \{\kappa_2\} \cup \mathfrak{J} & \text{otherwise} \end{cases}$
\end{center}
with $\mathfrak{J}$ as required in the proof of Proposition 3. One idea for applying Theorem 2 to $\models^*_\Delta$ may thus be that of modifying Definition 9 so as to have justifications defined, not from argument structures to argument structures, but from subsets of $\mathbb{S} \times \mathbb{B}$ to arguments, where $\mathbb{S}$ is the sets of the argument structures over $\mathscr{L}$ and $\mathbb{B}$ is the class of atomic bases. In this way, we may set
\begin{center}
    $\texttt{Ch} \langle \mathscr{D}, \mathfrak{B} \rangle = \langle \mathscr{D}, \mathfrak{J}^\mathfrak{B} \rangle$.
\end{center}
$\texttt{Ch}$ would thus behave like a choice function picking the right justifications set on each base, so $\langle \mathscr{D}, \texttt{Ch} \rangle$ is valid on every $\mathfrak{B}$ -- we remark that this strategy requires a number of changes in the formal apparatus of Section 2, which eventually lead to ToG where, roughly, justifications of rules are embedded into inference steps -- see \cite{piccolomini}.

Another alternative would be to put all the $\mathfrak{J}^\mathfrak{B}$-s into the justifications set for $\mathscr{D}$, so that $\langle \mathscr{D}, \bigcup_{\mathfrak{B} \in \mathbb{B}} \mathfrak{J}^\mathfrak{B} \rangle$ be valid on every $\mathfrak{B}$. Contrarily to the previous solution, this move does not require we touch Definition 9 at all.

A final option, which is similar to that where we take the union-set of the justifications set of each base, but which, similarly to the choice function strategy, involves a modification of Definition 9, amounts to following \cite{schroeder-heister2006}, namely, replacing, so to say, the justifications with their graph.\footnote{There are here some subtleties which I cannot deal with. In the approach inspired by \cite{schroeder-heister2006}, one normally allows for alternative justifications for one and the same argument structure -- a possibility I have left open in Definition 9. In general, one can show that from any justifications set one can extract what, in Definition 17 below, I call a r-sequence; if alternative justifications are permitted, the inverse holds too. Something similar obtains concerning the relationship between validity of Definition 12 and what I call SH-validity in Definition 20 below: if alternative justifications are not permitted, then the implication holds only for closed argument structures, otherwise one has an equivalence between the notions.} Let me give in this case more details than I have done so far. First of all, the modified Definition 9 now reads as follows.

\begin{definition}
    A \emph{reduction} is a pair $\langle \mathscr{D}^1, \mathscr{D}^2 \rangle$ where $\mathscr{D}^2$ has the same conclusion, and at most the same assumptions as $\mathscr{D}^1$. A \emph{r-system} is a set of reductions. A \emph{r-sequence} $\langle \mathscr{D}^1_1, \mathscr{D}^2_1 \rangle, ..., \langle \mathscr{D}^1_n, \mathscr{D}^2_n \rangle$ such that, for every $i \leq n$, $\mathscr{D}^2_{i - 1} = \mathscr{D}^1_i$, is said to be \emph{from} $\mathscr{D}^1_1$ \emph{to} $\mathscr{D}^2_n$.
\end{definition}
\noindent Let us modify accordingly Definitions 10, 11 and 12.
\begin{definition}
    Given a r-system $\Sigma$, $\mathscr{D}$ \emph{reduces to} $\mathscr{D}^*$ \emph{relative to} $\Sigma$, written $\mathscr{D} \leq^\Sigma \mathscr{D}^*$, iff $\Sigma$ contains a r-sequence from $\mathscr{D}$ to $\mathscr{D}^*$.
\end{definition}
\begin{definition}
    An \emph{argument} is a pair $\langle \mathscr{D}, \Sigma \rangle$.
\end{definition}
\begin{definition}
    $\langle \mathscr{D}, \Sigma \rangle$ is \emph{SH-valid} on $\mathfrak{B}$ iff the same conditions as in Definition 12 hold, with $\mathfrak{J}$ and $\leq^\mathfrak{J}$ replaced by $\Sigma$ and $\leq^\Sigma$ respectively.
\end{definition}
\noindent It is now easy to see that the following obtains.
\begin{proposition}
For every $\mathfrak{B}$, there is closed $\langle \mathscr{D}, \Sigma \rangle$ for $A \vee \neg A$ SH-valid on $\mathfrak{B}$.
\end{proposition}
\begin{proof}
    Basically the same as proof of Proposition 3, where the relevant reductions are $\langle \mathscr{D}, \kappa_1(\mathscr{D}) \rangle$ and $\langle \mathscr{D}, \kappa_2(\mathscr{D}) \rangle$, respectively.
\end{proof}
\noindent This shows that, for every $\mathfrak{B}$, there is a r-system $\Sigma_\mathfrak{B}$ such that $\langle \mathscr{D}, \Sigma_\mathfrak{B} \rangle$ is valid on $\mathfrak{B}$ -- where $\mathscr{D}$ is the axiom $A \vee \neg A$. Let us now give the following definition of logical validity.
\begin{definition}
    $\Gamma \models_\Delta^{SH} A$ iff, for some $\langle \mathscr{D}, \Sigma \rangle$ from $\Gamma$ to $A$ and every $\mathfrak{B}$, $\langle \mathscr{D}, \Sigma \rangle$ is SH-valid on $\mathfrak{B}$.
\end{definition}
\noindent Thus we have what follows.
\begin{proposition}
    $\models^{SH}_\Delta A \vee \neg A$.
\end{proposition}
\begin{proof}
    Take $\langle \mathscr{D}, \Sigma \rangle$ with $\mathscr{D}$ the axiom $A \vee \neg A$, and $\Sigma = \bigcup_{\mathfrak{B} \in \mathbb{B}} \Sigma_\mathfrak{B}$.
\end{proof}
\begin{theorem}
    There are $\Gamma$ and $A$ such that $\Gamma \models^{SH}_\Delta A$ and $\Gamma \not\vdash_{\texttt{IL}} A$.
\end{theorem}

There is however a clear sense in which neither $\texttt{Ch}$, nor $\bigcup_{\mathfrak{B} \in \mathbb{B}} \mathfrak{J}_\mathfrak{B}$ or $\bigcup_{\mathfrak{B} \in \mathbb{B}} \Sigma_\mathfrak{B}$ can be said to be or to consist of \emph{schematic} justifications, where by schematic we may provisionally mean that the justification can be expressed as a \emph{rule for rewriting argument structures}, similar to the reduction for removing disjunction detours that we have seen at the beginning of Section 2. $\texttt{Ch}$ is obviously not so, since it is not even defined on argument structures only, but on pairs of arguments structures \emph{and} bases. Even if it is in a sense schematic, it hence differs significantly in spirit from $\phi_\vee$.

On the other hand, $\bigcup_{\mathfrak{B} \in \mathbb{B}} \mathfrak{J}^\mathfrak{B}$ contains of course a schematic justification, i.e. $\kappa_1$. But the $\kappa_2$ we take on each base cannot be assumed to be schematic. This is a constant function which points the axiom $A \vee \neg A$ to a fixed closed argument structure for $A \vee \neg A$ obtained by introducing disjunction below a closed argument structure for $A$, where the latter is given. Clearly, we have no guarantee that the closed argument structure thereby obtained has a form that is invariant over all bases, and hence that $\kappa_2$ can be described as a rewriting scheme. Also, we have no information at all about $\mathfrak{J}$, i.e. the justifications sets for the closed argument structure for $A$ (if any). For essentially the same reasons, we may not be able to describe $\Sigma$ in proof of Proposition 6.

Under a suitable description of this notion of schematicity, we may thus modify $\models^*_\Delta$ in Definition 16, so as to obtain a further notion of logical consequence.
\begin{definition}
    $\Gamma \models_\Delta^s A$ iff, there is $\mathscr{D}$ from $\Gamma$ to $A$ such that, for some schematic $\mathfrak{J}$ and every $\mathfrak{B}$, $\langle \mathscr{D}, \mathfrak{J} \rangle$ is valid on $\mathfrak{B}$.
\end{definition}
\noindent This somewhat stricter reading may cope with constructivist desiderata. The validity of excluded middle may be no longer provable, even with classical meta-logic.

\section{Conclusion}

Prawitz's original aim was that of developing a constructive semantics which accounted for $\texttt{IL}$ \cite{prawitz73}. Investigations into proof-theoretic semantics have instead shown that $\texttt{IL}$ is in general \emph{not} complete, or that it is so only on condition that atomic bases undergo certain constraints -- see \cite{piecha} for an overview, while more recent findings are \cite{piechaschroeder-heister, schroeder-heisterolf, stafford1, stafford2}. However, these incompleteness (or restricted completeness) results are normally referred to a Base Semantics approach, where justifications and proof-structures are dropped out. In some cases, such results can be adapted to the SVA approach, but this may require understanding reductions in a non-schematic way, say, as choice-functions defined also on bases, or as non-recursive sets of justifications or of reductions sequences in Schroeder-Heister's sense. If we require more schematicity, the adaptation may fail, e.g., even with classical logic in the meta-language, we may no longer be able to prove the validity of excluded middle in the object-language. This would speak in favour of Prawitz's original semantic project being constructive in two strictly interrelated ways, namely, not just because SVA is based on a notion of valid argument (i.e. of proof), but \emph{also} because this notion is given in terms of reductions which amount to concrete rewriting rules for proof-structures.

Of course, the question now becomes how schematicity should be more precisely defined. I will conclude by hinting at two possible strategies for dealing with this topic.

First, one may consider the possibility of relaxing the notion of open validity, i.e. an open valid argument $\langle \mathscr{D}, \mathfrak{J} \rangle$ is such that $\mathscr{D}$ reduces to canonical form relative to $\mathfrak{J}$ even when its assumptions \emph{are not} replaced by closed structures which be valid relative to expansions of $\mathfrak{J}$. Observe that this \emph{does not hold} in SVA. For example, take the base given by the rules
\begin{prooftree}
    \AxiomC{}
    \UnaryInfC{$p$}
    \AxiomC{$p$}
    \UnaryInfC{$q$}
    \noLine
    \BinaryInfC{}
\end{prooftree}
and consider the rule
\begin{prooftree}
    \AxiomC{$r$}
    \UnaryInfC{$q \vee s$}
\end{prooftree}
justified by a $\phi$ such that
\begin{prooftree}
    \AxiomC{}
    \UnaryInfC{$p$}
    \UnaryInfC{$r$}
    \UnaryInfC{$q \vee s$}
    \AxiomC{$\Longrightarrow$}
    \noLine
    \UnaryInfC{}
    \AxiomC{}
    \UnaryInfC{$p$}
    \UnaryInfC{$q$}
    \RightLabel{($\vee_I$)}
    \UnaryInfC{$q \vee s$}
    \noLine
    \TrinaryInfC{}
\end{prooftree}
Then, $\langle \mathscr{D}, \phi \rangle$ is valid on our base, where $\mathscr{D}$ is
\begin{prooftree}
    \AxiomC{$p$}
    \UnaryInfC{$r$}
    \UnaryInfC{$q \vee s$}
\end{prooftree}
It may nonetheless be the case that open structures whose assumptions have some \emph{peculiar feature} may be computable to canonical form relative to some base-independent justifications, namely, justifications which be schematic in the broad sense understood here. A solution of this kind may come from \emph{Pezlar's selector} \cite{pezlar} for the \emph{Split rule}
\begin{prooftree}
    \AxiomC{$A \rightarrow B \vee C$}
    \UnaryInfC{$(A \rightarrow B) \vee (A \rightarrow C)$}
\end{prooftree}
with $A$ Harrop formula. Pezlar's selector is based on the idea that open proofs of the implicational premise involve sufficient computational content for them to reduce to canonical form under trivial \emph{ad hoc} proofs of the antecedent.

Alternatively, one may put a sort of linearity constraint on justifications, like the one put by Definition 9 on applications of justifications to instances of open structures. We may require the same for the structures which the justification is defined on. If we go back to the reduction for disjunction detours, we see that it consists of a \emph{constant} part, i.e. the order of the formulas on the left- and right-hand side of the arrow, and of a \emph{variable} part, given by variables for argument structures. Let us express this as a linear operation $\vee E$ on typed proof-objects as happens in ToG, and let us indicate with $\texttt{inj}$ a function for obtaining a proof of $A_1 \vee A_2$ from one of either $A_1$ or $A_2$. Then, the rule enjoys the following linearity condition: for any $f^1_1, f^1_2 : (A_1)B$, $f^2_1, f^2_2 : (A_2)B$ and $x^i_1, x^1_2 : A_i$ ($i = 1, 2$),
\begin{center}
    $\vee E(\texttt{inj}(x^i_1)[x^i_2/x^i_1], f^1_1(y_1)[f^1_2/f^1_1], f^2_1(y_2)[f^2_2/f^2_1]) =$
\end{center}

\begin{center}
    $= (\vee E(\texttt{inj}(x^i_1), f^1_1(y_1), f^2_1(y_2)))[x^i_2, f^1_2,f^2_2/x^i_1, f^1_1, f^2_1] =$
\end{center}

\begin{center}
    $= f^i_1(x^i_1)[x^i_2, f^1_2, f^2_2/x^i_1, f^1_1, f^2_1] = f^i_1(x^i_1)[x^i_2, f^i_2/x^i_1, f^i_1]$.
\end{center}
where $y_1 : A_1$, $y_2 : A_2$, and both are bound by $\vee E$. This may become a definitory condition for justifications, that is, we only allow for justifications $\phi$ such that, for every replacement $\sigma$ of proof-variables in the argument $x$ of $\phi$,
\begin{center}
    $\phi(x^\sigma) = (\phi(x))^\sigma$.
\end{center}

\end{document}